\documentclass{llncs}%
\usepackage{amsfonts}
\usepackage{amssymb}
\usepackage[nosumlimits,nointlimits,nonamelimits]{amsmath}
\usepackage[onehalfspacing]{setspace}
\usepackage{footmisc}
\usepackage{graphicx}
\usepackage{comment}
\input{epsf.sty}
\input xypic
\xyoption{all}
\LaTeXdiagrams
\xyoption{2cell}
\UseAllTwocells
\providecommand{\U}[1]{\protect\rule{.1in}{.1in}}
\UseAllTwocells

\newtheorem{teo}{Theorem}

\newcommand{\f}[2]{\frac{\displaystyle #1}{\displaystyle #2}}
\def\sq{\sqrt}
\def\sq2{\sqrt{2}}
\def\sq12{\sq{12}}
\def\dsq2{\f{1}{\sqrt{2}}}
\def\be{\begin{equation}}
\def\ee{\end{equation}}
\def\q={\quad = \quad}
\def\lra{\longrightarrow}

\def\A{{\mathbb A}}
\def\B{{\mathbb B}}
\def\C{{\mathbb C}}
\def\D{{\mathbb D}}
\def\Cat{{\mathbb Cat}}

\def\N{{\mathbb N}}

\def\X{{\mathbb X}}

\newcounter{examnum}[section]
\newcounter{remarnum}[section]
\begin{document}
\title{Kolmogorov Complexity of Categories}
\author{Noson S. Yanofsky}
\institute{Department of Computer and Information Science,
Brooklyn College, CUNY, Brooklyn, N.Y. 11210.\\ and the Computer Science
Department of the Graduate Center, CUNY,\\ New York, N.Y. 10016. \\
\email{noson@sci.brooklyn.cuny.edu}\thanks{A while back, I showed some of these ideas to Samson Abramsky and he was, as always, full of encouragement and great ideas. I am very grateful to him for all his help over the years. I would like to acknowledge the help and advice of  Michael Barr, Marta Bunge, James Cox, Joey Hirsh, Florian Lengyel, Dustin Mulcahey, Philipp Rothmaler, and Louis Thral. I want to thank Shayna Leah Hershfeld for many enlightening conversations about polymorphism and type theory. Support for this project was provided by a PSC-CUNY Award, jointly funded by The Professional Staff Congress and The City University of New York.}}

\maketitle
\begin{abstract}
 Kolmogorov complexity theory is used to tell what the  
algorithmic informational content of a string is. It is defined as the length of the
shortest program that describes the string. We present a programming
language that can be used to describe categories, functors, and
natural transformations. With this in hand, we define the
informational content of these categorical structures as the shortest
program that describes such structures. Some basic consequences of our definition are
presented including the fact that equivalent categories have equal Kolmogorov complexity. We also prove different theorems about what can and cannot be described by our programming language. 

\vspace{.1in}


\vspace{.1 in}

\noindent {\bf Dedicated to Samson Abramsky in honor of his 60th Birthday}
\end{abstract}

\section{Introduction}
Kolmogorov complexity is a part of theoretical computer science that was pioneered in the early 1960's by Andrey Kolmogorov, Ray Solomonoff, and Gregory Chaitin. For reasons ranging from probability theory, to machine learning, and computational complexity theory, these three researchers gave a universal definition of what it means for a string of symbols to be simple or complex.

Consider the following three strings:

1. 00000000000000000000000000000000000000000000000

2. 11011101111101111111011111111111011111111111110

3. 01010010110110101011011101111001100000111111010

\noindent All three consists of 0’s and 1’s and are of length 45.  It should be noticed that if you flipped a coin 45 times the chances of getting any of these three sequences are equal. That is, the chances for each of the strings occurring is $1/2^{45}.$ In effect, this shows a failure of classical probability theory in measuring the contents of a string. Whereas you would not be shocked to see a sequence of coins produce string 3, the other two strings would be surprising. The difference between these strings can be seen by looking at short programs that can describe them: 

1. \verb"Print 45 0's".

2. \verb"Print the first 6 primes."

3. \verb"Print `01010010110110101011011101111001100000111111010'."

\noindent The shorter the program, the less informational content of the string. In contrast, if only a long program can describe the string, then the string has more content.
If no short program can describe a string, then it is ``incompressible" or
``random."

In classical Kolmogorov complexity, rather than talking about programs, one talks about Turing machines. For a string $s$, the   the Kolmogorov complexity, $K(s)$, is defined as the size of the smallest Turing machine that starts with an empty tape and outputs $s$.
Formally, let $U$ be a universal Turing machine, then
$K(s)= min \{ |p| : U(p,\lambda)=s \}. $
We will also need {\it relative} Kolmogorov complexity: let $s$ and $t$
be two strings, then $K(s|t)$ is the size of the smallest Turing machine
that starts with $t$ on the tape and outputs s.
Formally,
$K(s|t)= min \{ |p| : U(p,t)=s \}. $ If $K(s)>|s|$ then $s$ is ``incompressible'' or ``random''.

This notion of Kolmogrov complexity is used in many different areas of theoretical computer science. It gives an objective measure of how complicated strings are. 
It is our goal to extend these ideas to many other areas of mathematics, computer science and physics by formulating a notion of Kolmogorov complexity for category theory which is used in all these diverse areas. In order to measure how complicated categories, functors, and natural transformations are, we need a programing language that will describe these categorical structures. In honor of Sammy Eilenberg, one of the founders of category theory who also had a deep interest in computer science, we call this programming language ``Sammy.'' This language will have variables that can hold categories, functors and natural transformations. The operations of the language will perform common constructs that people use to formulate different structures. Each line of the program could have a label that will be used with "If-Then" statements to control the execution of the program.

Notice that numbers, strings, trees, graphs, arrays, and other typical data types are not mentioned in our programming language. This was done on purpose. The other data types can be derived from the categorical structures. Categories and algorithms are more ``primitive'' than numbers, strings, etc. 

This is not the first time a programing language has been formulated to describe categorical structure.  An important example is in {\it Computational Category Theory} by Rydeheard and Burstall \cite{RandB}.  Tatsuya Hagino's thesis \cite{Hagino} is another example. These languages are, however, different from Sammy. Their programming languages are made to be implemented and to get computers to actually calculate with categories. In contrast, there is no intention of implementing Sammy. Our goal is simply to compare different structures by comparing the length of their descriptions. In fact, we will not even write many formal Sammy programs. This is similar to the fact that no one actually ever formally writes the instructions for a Turing machine.  

With Sammy, we will talk about the Kolmogorov complexity of categorical structures. We discuss when  one structure is more complicated than another. We will also talk about compressibility and randomness. Along these lines, here is a simple example of the type of ideas we will meet. Consider $\N$, the totally ordered category of natural numbers $\xymatrix{0\ar[r]&1\ar[r]&2\ar[r]&\cdots}$, and ${\bf \overline{2}}$, the category with two objects and a single isomorphism between them  $\xymatrix{0\ar[r]^{\sim}&1}$. A functor $F:\N \lra {\bf \overline{2}}$ corresponds to an infinite sequence of zeros and ones. The category of all such functors ${\bf \overline{2}}^{\N}$ is essentially to the real numbers and has uncountably many elements. How many of these functors can be mathematically described? There are only countably many computer programs that describe such functors. This means that the vast majority of functors $\N \lra {\bf \overline{2}}$ cannot be described by any program and are essentially random.

Not every categorical structure can be described with our programing language. Categorical structures that can be described by Sammy will be called ``constructible.'' For example, I do not know how to start from nothing and make the category of smooth manifolds. However it is probably possible to start from the category of topological spaces and get the category of smooth manifolds. This brings us to the notion of relative Kolmogorov complexity. We will be interested in how long does a program have to be in order to construct a categorical structure {\it given some categorical structures.}

The fact that certain structures are not constructable with Sammy brings in the whole area of computability theory. There are limitations to what Sammy can perform. 
Usual self-referential limitations are based on variations of the liar paradox (``This statement is false'') such as
G\"odel (``This statement is unprovable'') or Turing (``This program will output the wrong answer when asked if it will halt or go into an infinite loop'') (see \cite{OLR} for a comprehensive survey of such limitations.) In contrast, the limitations of Kolmogorov complexity 
are based on the Berry Paradox: consider the number described by ``The least number that needs more than fifteen
words to describe it.'' This sentence has twelve words. That is, there is a description of a number that is shorter than it is supposed to be. One such limitation within classical Kolmogorov complexity\cite{Vitanyi} is:

\begin{teo} $K:\mbox{Strings}\lra \mathbb N$ is not a computable function.
\end{teo}

\noindent We will show that there are similar limitations for our Kolmogorov complexity theory.

Section 2  introduces Sammy. That section also describes several ``library functions'' or ``macros'' in Sammy which will be helpful in the rest of the paper. Section 3 is the heart of the paper where we define and prove many of the central theorems about our complexity measure. 
Section 4 is a discussion of computability and non-computability with the Sammy language. The paper concludes with some possible ways this work will progress in the future.

\section{A Programing Language for Categories}

In order to describe categorical structures, we need a programing language. This language will be called ``Sammy''.
The language will consist of typical operations that are used to describe/create different
categories, functors and natural transformations. Programs will be lists of statements that set variables to different values. The variables could be categories, functors, or natural transformations. Since categories are special types of functors, and functors are special types of natural transformations (that is, natural transformations are the deepest type), we might state everything in terms of natural transformations. But that would make the programs needlessly complex. Rather, for the sake of simplicity, we will be ambiguous about the types of our statements (that is, our operations/functions will be polymorphic.)
As we have absolutely no intention of implementing Sammy, we can be vague about certain issues.

We begin with constants. There is $\bf 0$, the empty category, $\bf 1$,  the category with one object and one morphism, and $\bf 2$, the category $0\lra 1$ with two objects and one nontrivial morphism. We will also need the constant category \verb"Cat" which corresponds to the category of all small categories. There are also several constant functors: $s:\bf 1 \lra \bf 2$ and $t:\bf 1 \lra \bf 2$ that picks out the source and target of the nontrivial morphism in $\bf 2$. There are the unique morphisms $!:\bf \bf 0 \lra \bf 1$, $!:\bf \bf 0 \lra \bf 2$,
$!:\bf 0 \lra \verb"Cat"$, $!: \verb"Cat" \lra \bf 1$, and $!: \bf 2 \lra \bf 1$. There are also identity functors and natural transformations.  

There are several operations that take a single input. For a functor $F:\A \lra \B$ if we set 
$\C=\verb"Source"(F:\A \lra \B)$ then $\C=\A$. That is, \verb"Source" takes a functor and outputs the category that is the source of the functor. There is a similar operation $\C=\verb"Target"(F:\A \lra \B)$. For a given category $\A$, the operation
$F=\verb"Ident"(\A)$ makes $F=Id_\A$. For a category $\A$, if we let $\C=\verb"Op"(\A)$ then $\C=\A^{op}$. The \verb"Op" operation also acts on functors.

We will at times have to talk about an actual object and morphism in the category. So for example, a functor $F:\bf 1 \lra \C$ ``picks'' an object $c$ in $\C$ and a functor $F: \bf 2 \lra \C$ ``picks'' a morphism $f:c \lra c'$. Going the other way, an object $c$ in $\C$ ``determines'' a functor $F_c: \bf 1 \lra \C$ and similarly for a morphism in $\C$. We write this in Sammy as $c=\verb"Pick"(F:\bf 1 \lra \C )$ and $F_c=\verb"Determine"(c)$.   

For natural transformations of the appropriate source and target there is a horizontal composition and vertical composition written as  $\alpha=\verb"Hcomp"(\beta, \gamma)$ and  $\alpha=\verb"Vcomp"(\beta, \gamma)$. Regular composition of functors is simply a special case of horizontal composition. For categories $\A$ and $\B$, we will have $\C = \verb"Pow"(\A, \B)$ be the category of all functors and natural transformations from $\A$ to $\B$. 

Probably the most important operations are the Kan extensions. For functors $G:\A \lra \B$ and $F:\A \lra \C$, a right Kan extension of $F$ along $G$ is a pair $(R, \alpha)=\verb"KanEx"(G, F)$
where $R:\B \lra \C$ and $\alpha:R\circ G \lra F$. A Kan extension induces another functor that is unique. For every $H:\B \lra \C$ and $\beta: H \circ G \lra F$ there is a unique $\gamma=\verb"KanInd"(F,G; H, \beta)$
where $\gamma:H \lra R$ and satisfies $\alpha \cdot \gamma_G = \beta$.  
Using Kan extensions one can derive, products, coproducts, pushouts, pullbacks, equalizers, coequalizers, (and constructible) limits, colimits, ends, coends, etc. It is a well-known fact that if $G:\A \lra \B$ is a right adjoint (left adjoint, equivalence, isomorphism), then its left adjoint (right adjoint, quasi-inverse, inverse) $G^*:\B \lra \A$ can be found as a simple Kan extension of the identity $Id_\A$ along $G$, that it, $G^*=\verb"KanEx"(G, Id_{\A})$.

For ``bootstrapping'' purposes we will need an operation that takes two categories and gives their coproduct and their induced maps. This will help us create categories like $\bf 1 \sqcup \bf 1$ which will be needed for our Kan extensions to describe products and coproducts; and $\bf 2 \sqcup \bf 2$ which will be needed to describe equalizers and coequalizers.   

There is a dual notion of a Kan Lifting. For functors $F:\A \lra \B$ and $G:\C \lra \B$ a Kan lifting of $F$ along $G$ is a pair 
$(R,\alpha) = \verb"KanLif"(G, F)$
where $R:\A \lra \C$ that satisfies a universal property which can easily be written down.

Since Kan extensions and Kan liftings are only defined up to a unique isomorphism, we might ask what is the output of the function $\verb"KanEx"(G, F)$? We do not care. The computer decides which of the many possible outputs it will output. It is irrelevant from the categorical perspective. This is similar to a real programing language when we do not know how something is stored or how a function is calculated. The user is ambivalent as to how the computer does certain actions. 
We are also well-aware that the Kan extensions and Kan liftings might not exist. In that case, the program will not go on.

There is one more operation that needs to be discussed. Let $\C$ be a category. $\C^{\bf 2}$ and $\C^{\bf 1}$ are the categories of arrows and objects of $\C$. The maps $s:\bf 1 \lra \bf 2$ and $t:\bf 1 \lra \bf 2$ induce (using the \verb"Pow" operation on functors) maps $\C^s:\C^{\bf 2} \lra \C^{\bf 1}$ and $\C^t:\C^{\bf 2} \lra \C^{\bf 1}$. The pullback of these two maps, $\C^{\bf 2} \times_{\C^{\bf 1}} \C^{\bf 2}$ is the composable arrows in the category. The important part of the  information about the category is the composability map $\circ:  (\C^{\bf 2} \times_{\C^{\bf 1}} \C^{\bf 2}) \lra \C^{\bf 2}$. This map will help us get into the nitty-gritty of how a category is defined. So we have the following operation: for a category $\C$, the operation $F=\verb"Composable"(\C)$ gives us the $\circ$ map.  

We would like some control of how the Sammy program will execute. We do this with a conditional branch statement: \verb"If " $\alpha_1==\alpha_2$ \verb"Goto L" where $\alpha_1$ and $\alpha_2$ are natural transformations and \verb"L" is a label of some program line. With such a conditional branch, we can get all the usual logical operations: AND, NOT, etc. We can also get the unconditional branch \verb"Goto L".  

There are a number of remarks that need to be made about Sammy:

This might not be the best language for our purposes.  Certain operations can be derived from other operations and hence a smaller more compact language is possible. For example, the \verb"Target" operation can be derived from the \verb"Source" and \verb"Op" operations. Bear in mind that our goal is to count the number of operations up to a coefficient. So we need not be exact. If one operation can be replaced by a constant number of other operations, nothing is lost. 

This language can not describe all constructions. (We shall see later.) What can be
done with this language will be called ``constructible.''
It is interesting to look at what type of categories can be described by
this programming language with no other input.

There is a need for a Church-Turing type thesis. The classic Church-Turing thesis says that whatever can be computed, can be computed by a Turing machine. We need such a thesis that says that whatever can be constructed by categorical means, can be constructed using the Sammy programing language. Alas, this is a thesis and not a theorem because we cannot characterize what can be constructed by categorical means. We will see that there are certain constructions that cannot be performed by Sammy. However, we believe that no programming language can make those constructions. 

With classical Kolmogorov complexity, there is much discussion about ``self-delimiting'' programs. This will not be an issue here. We can easily tell when a Sammy program begins and when it ends. 

With Sammy in hand, we introduce some library functions or macros that will be used in the future:

The coequalizer $ \xymatrix{ \mathbf 1 \ar@<-.1in>[r]_s \ar@<.1in>[r]^s& \mathbf 2 \sqcup \bf 2 }$ gives the category $\xymatrix{ \ast &\ast \ar[r] \ar[l]&\ast}$ which can be put in a Kan extension and give us pushouts and pullbacks. 
We can make many similar constructions.

For functors $L:\A \lra \C$ and $R:\B \lra \C$ we can construct the comma categories as the following pullbacks:
$$
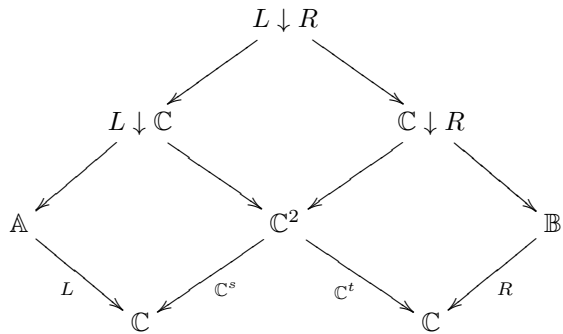
\xymatrix{&& L \downarrow R \ar[ld]\ar[rd] \\
 & L \downarrow \C \ar[rd]\ar[ld] &   &\C \downarrow R \ar[rd]\ar[ld]&\\
\A\ar[rd]_L&&\C^2\ar[rd]_{\C^t}\ar[ld]^{\C^s}&&\B\ar[ld]^R \\
 &\C && \C
}$$
Special instances of comma categories are slice categories and coslice categories.

The coequalizer
$ \xymatrix{ \mathbf 1 \ar@<-.1in>[r]_t \ar@<.1in>[r]^s& \mathbf 2 \ar[r]^{\rho}& \omega}$
gives the (infinite) natural numbers as a monoid. $\N= \mathbb \omega^{\mathbf 2}$ gives the totally ordered category of natural numbers. The successor function is defined as follows:
$$\xymatrix{r:\omega\ar[r]^{\sim}&\omega\times \bf 1 \ar[r]^{Id \times s}& \omega \times \bf 2 \ar[r]^{Id \times \rho}& \omega \times \omega \ar[r]^{\circ}& \omega .}$$
That is, take any $n\in \omega$ and associate it with the nontrivial morphism in $\bf 2$. This becomes the $+1$ member of $\omega$. Then compose $n$ with $+1$. Now take this map $r$ and look at $s=r^{\bf 2}: \N=\omega^{\bf 2} \lra \omega^{\bf 2}=\N.$ This is the successor map. 

We construct the category with two objects and a unique isomorphism between them. First make a category with two distinct copies of $\bf 2$. By keeping track of the inclusion maps, we have an induced $F$ and $G$
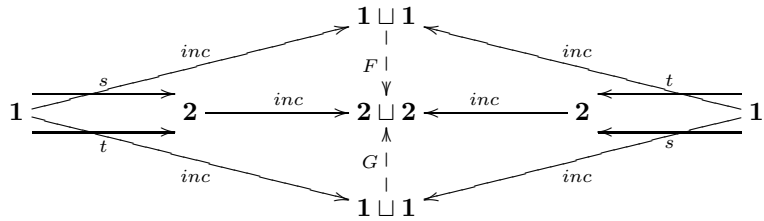
$$ \xymatrix{
&&&& \mathbf 1\sqcup \mathbf 1 \ar@{-->}[d]_F\\ 
\mathbf 1 \ar@<-.1in>[rr]_t \ar@<.1in>[rr]^s\ar[drrrr]_{inc}\ar[urrrr]^{inc}&& \mathbf 2 \ar[rr]^{inc}&& \mathbf 2\sqcup\mathbf 2
&&\mathbf 2 \ar[ll]_{inc} && \mathbf 1 \ar@<-.1in>[ll]_t \ar@<.1in>[ll]^s\ar[dllll]^{inc}\ar[ullll]_{inc}\\ 
&&&& \mathbf 1\sqcup \mathbf 1 \ar@{-->}[u]^G}$$
Now use these induced maps in a coequalizer to form the desired category.  The figure on the right is helpful.
$$ 
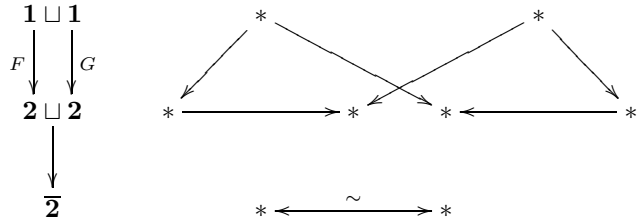
\xymatrix{\mathbf 1 \sqcup \mathbf 1\ar@<-.1in>[d]_F \ar@<+.1in>[d]^G& & \ast\ar[dl]\ar[drr] & & &\ast\ar[dr]\ar[dll]\\ 
\mathbf 2 \sqcup \mathbf 2\ar[d]& \ast \ar[rr] && \ast& \ast&& \ast \ar[ll]\\
\overline{\mathbf 2}&&\ast \ar[rr]^{\sim}&&\ast \ar[ll]
}$$

\section{Kolmogorov Complexity of Categories}

For a category $\mathbb{C}$ (or a functor, or a natural transformation) we define $K_{Sammy}(\mathbb{C})$ to be the number of operations in the smallest Sammy
program that describes $\mathbb{C}$. For relative Kolmogorov complexity, letting $$\Gamma=\{\mathbb{C}_1,\mathbb{C}_2, \ldots ,\mathbb{C}_l,F_1,F_2, \ldots, F_m, \mu_1, \mu_2, \ldots ,\mu_n\},$$
or $\Gamma$ as a sub2-category of $\mathbf{Cat}$ then $K_{Sammy}(\mathbb{C}|\Gamma)$ is the number of operations in the smallest Sammy program that describes $\mathbb{C}$ given $\Gamma$ as input. We shorten $K_{Sammy}$ to $K$ when no confusion will arise.

If there is a finite number of operations so that one can go from one categorical structure to another and vice versa, we say that the Kolmogorov complexity of these categorical structures are approximately the same. In detail, if there exists a $c$ such that for all appropriate categorical structures, $\X$, one can change $\X$ to $\X'$ and vice versa in $c$ Sammy operations, that is $|K(\X) - K(\X')|\leq c$, then we write $K(\X) \approx K(\X')$. As an example, notice that only one Sammy operation is needed to go from category $\A$ to functor $Id_{\A}$ and vice versa. Hence $K(\A) \approx K(Id_{\A})$. 

There is a need for something called an {\it invariance theorem}. This basically says that the Kolmogorov complexity does not depend on the programing language that is used to describe the objects. Imagine that you do not like the Sammy programing language to describe categorical structures and you decide to invent your own. Perhaps you call it ``Saunders'' (after the other founder of category theory, Saunders Mac Lane.) Then since presumably both languages can program any constructable categorical structure, they can each program the other's operations. That means there exist compilers that can translate Sammy programs into Saunders programs and there are compilers that can translate Saunders programs into Sammy programs. From this, we can prove the following theorem: 
There exists a constant $c$ such that for all categorical structures $\bf X$ we have $|K_{Sammy}({\bf X}) - K_{Saunders}({\bf X})|\leq c$.

Rather than list all the results we have for $K$, let us examine some paradigmatic theorems:

\begin{teo} There exists a constant $c_{pair}$ such that for all $\mathbb C$ and $\mathbb D$ we have
$K({\mathbb C \times \mathbb D}) \le K(\mathbb C)+K(\mathbb D|\mathbb C)+c_{pair}$. \end{teo}

This essentially says that there is a simple way of taking two categories and forming their product. There is no new information added. But lets look more carefully at what the theorem say. It says that to form $\C \times \D$ one can form $\C$ and then form $\D$ (but you might use some information that you already have since you already formed $\C$) and then do a few lines of Sammy to get their product. The reason for the inequality is because there might be an easier way. For example $\bf 0 \times \D$ can be formed in a constant amount of operations: it is $\bf 0$. There is also a similar theorem with $\C$ and $\D$ swapped on the right side of the inequality. 

\begin{teo} There exists a constant $c_{double}$ such that for all $\mathbb C$ we have $K({\mathbb C \times \mathbb C}) \le K(\mathbb C)+c_{double}$ .\end{teo}
That is, there is a simple way to double a category and no new information is there.

\begin{teo} There exists a constant $c_{target}$ such that for all $F:\A \lra \B$ we have $K(\B) \le K(F:\A \lra \B)+c_{target}.$
\end{teo}

This means that one way to describe $\B$ is to first find a program for a functor $F:\A \lra \B$ and then use the \verb"Target" operation to get $\B$. The inequality comes from the fact that there might be shorter programs to describe $\B$.
There are similar such theorems for the source of a functor, for natural transformations, for identity functors, etc. 

We state the following theorem about composition in terms of natural transformations for generality.
\begin{teo} There exists a constant $c_{compos}$ such that for any three natural transformations 
$\alpha : F \lra G$, $\beta : F \lra H$, and $\gamma:G \lra H$ such that $\beta= \gamma \circ \alpha$ we have 
$$K(\beta) \leq K(\alpha )+K(\gamma|\alpha) + c_{compos}. $$
\end{teo}

When $\gamma$ is the unique natural transformation that satisfies this triangle (e.g. when $\alpha$ is mono) then the inequality in the above theorem becomes an equality. 
 
The theorem for Kan extensions is similar. 
\begin{teo}
There exists a constant $c_{Kan}$ such that for all $G:\mathbb A \lra \mathbb B$ and
$F:\mathbb A \lra \mathbb C$ if $(Lan_G(F), \alpha)$ is the left Kan extension, than 
$$ K((Lan_G(F), \alpha)) \le K(F) + K(G|F)+ c_{Kan}$$
or for relative Kolmogorov complexity 
$$K((Lan_G(F),\alpha)|\Gamma)  \le K(F|\Gamma)+ K(G|\Gamma,F)+ c_{Kan}.$$
\end{teo}

As a special case, if $G:\A \lra \B$ is a right adjoint (left adjoint, equivalence, or isomorphism), then the Kan extension along $G$ of the $Id_{\A}$  is the left adjoint (right adjoint, quasi-inverse, inverse) $G^*:\B \lra \A$. Since it is easy to go from one to the other, we have that $K(G)\approx K(G^*)$. Notice that for an arbitrary adjunction, this does not mean that $K(\A) \approx K(\B)$ (we shall see that it is true for an equivalence). Nor does there seem to be any hard-and-fast rule that says something like a left adjoint goes from something with a low Kolmogorov complexity to a high Kolmogorov complexity. It is easy to find counterexamples to such ideas.


If $G:\A \lra \B$ and $F:\A \lra \C$ are functors, $R:\B \lra \C$ is a right Kan extension,  $H:\B \lra \C$, 
and $\beta: H\circ G \lra F$ then for the unique induced $\gamma:H \lra R$, we have that $K(\gamma) \approx K(\beta)$. The reason for this is that you can go from one to the other using composition and the \verb"KanInd" operation. A simple example of this is product:
$$\xymatrix{ &H\ar[rd]^{\beta_1}\ar[ld]_{\beta_0} \ar[d]^{! \gamma }\\
F_0 & F_0 \times F_1\ar[r]_{\alpha_1}\ar[l]^{\alpha_0} & F_1
}$$
It is easy to see that the information in $\gamma$ is exactly the information in the $\beta$s. It is easy to derive one from the other.



Our work would be in vain if the measure we described was not an invariant of categorical structure. We have the following important theorem.


\begin{teo} If categories $\A$ and $\B$ are equivalent, then $K_{Sammy}(\A)\approx K_{Sammy}(\B).$
\end{teo}

\noindent{\bf Proof.} The intuition behind the theorem is that Sammy cannot distinguish categorical structures that are isomorphic. Say the equivalence is given by the functor $G:\A \lra \B$. From $G$ its easily constructed quasi-inverse is $G^*: \B \lra \A.$  We then have that $K(G) \approx K(G^*)$. We also get that $K(G \circ G^*) \approx K(G^* \circ G)$. If $\alpha: Id_{\A} \lra GG^*$ is the isomorphic unit of the equivalence given by the Kan extension, then $\alpha^{-1}:GG^* \lra Id_{\A}$ is easily constructed (we are assuming that Kan extensions work on natural transformations). Since $\alpha^{-1} \circ \alpha = id_{Id}$ we get that $K(\alpha^{-1}) \approx K(Id_{\A})$ . 
We then have
$$K(\A) \approx K(Id_{\A})  \approx K(GG^*)  \approx K(G^*G) \approx K(Id_\B)\approx K(\B).$$
QED.

 There are some important consequences of this theorem. One can easily construct the skeletal category as the coequalizer
$ \xymatrix{ \mathbb C^{\overline{2}} \ar@<-.1in>[r]_t \ar@<.1in>[r]^s& \mathbb C \ar[r]^\cong& \mathbb C_{skeletal}}.$
This gives us $K(\mathbb C)\approx K(\mathbb C_{skeletal})$.

In a future paper \cite{YanoFut} we will discuss algebraic theories, monads, Morita equivalence and other algebraic notions from the Kolmogorov complexity perspective.

\section{Computability and Non-Computability with Sammy}
\noindent There might be a need to deal with finite numbers. We shall let the number $n$ correspond a triple $({\bf n}, P_b, P_e)$ where $\bf n$ is the totally ordered category with $n$ elements (keep in mind:  $\xymatrix{0 \ar[r]&1\ar[r]&\cdots \ar[r]&n-2\ar[r]&n-1}$), $P_b:\bf 1 \lra \bf n$ is a functor that points to the beginning of the category (the initial object), and $P_e:\bf 1 \lra \bf n$ is a functor that points to the end of the category (the terminal object.) Basic operations with such numbers are easy to describe. For example, we can connect $({\bf n}, P_b, P_e)$ and $({\bf m}, P'_b, P'_e)$ to get $({\bf n+m-1}, P_b, P'_e)$ with the coequalizer: $ \xymatrix{ \mathbf 1 \ar@<-.1in>[r]_{P_e} \ar@<.1in>[r]^{P'_b}& \mathbf n \sqcup \bf m \ar[r]& ({\bf n+m-1})}.$ (In truth, natural numbers can simply be given as functors $\bf 1 \lra \N$. We can manipulate numbers by manipulating such functors. While this is simple and economical, there is a certain appeal to doing it the way we did. Many prefer to think of their numbers as ``things'' and not just pointers to amounts.) 

All the finite totally ordered sets should be considered subcategories of $\N$ and, as such, inherit a partial successor function. Before applying this successor function we must check to make sure that the pointer is not at the $P_e$ position. 

A totally ordered category with $n$ elements can be constructed in $O(log_2 n)$ number of Sammy statements. Basically, the idea is that one can look at the binary representation of $n$ and write a program based on that. For example 
727 in binary is 1011010111. We can express this number as 
$$(((((((((1 \times 2+ 0)  \times 2+ 1) \times 2+ 1) \times 2+ 0) \times 2+ 1) \times 2+ 0)\times 2+1)\times 2+1)\times 2+1). $$
Similarly when making our totally ordered category, we can either (a) double the length of the category by connecting one copy of itself to itself, or (b) double itself and add one, depending on the bit at that position.  
This proves that
$K({\bf n}) \leq O(log_2 n)$
which is similar to the classical case.

Notice that the above algorithm did not have any input. In contrast, we can look at a program that loops through input, reads the bit and performs either (a) or (b). This input will be given as a functor from ${\bf log_2n}$ to $\overline{\mathbf 2}$. The program moves a pointer forward on ${\bf log_2n}$. There will be a conditional branch to see if the pointer is equal to $P_e$. While this might be a long program, it does not depend on the size of the input. We have thus proved that  
$$K({\bf n}\quad |\quad (F:{\bf log_2n} \lra \overline{\mathbf 2})) =O(1)$$
where $F$ describes $n$ in binary.


Considering numbers as such triples, we have the following theorem:
\begin{teo}  Any partially computable function of natural numbers can be computed with Sammy.
\end{teo}

\noindent{\bf Proof.} We prove that Sammy can perform the initial functions, recursion, composition, and the $\mu$-minimization operator. The zero function is achieved by simply setting $P_e=P_b$. The successor of $\bf n$ is achieved by simply composing with $\bf 2$. The projection function is simply a Sammy program that accepts $n$ inputs and outputs one of the inputs. 
Recursion can be done by iteration: we loop through a number until a pointer reaches $P_e$. Composition is simply composition of Sammy programs. $\mu$-minimization is done by doing a loop along $\N$ the ordered category of {\it all} natural numbers. QED.


What about complexity theory? In \cite{YanoFut} it is shown that categories and functors can mimic a Turing machine. For every rule of a Turing machine there is a set amount of steps of a Sammy program. Hence our programming language can do whatever a Turing machines can do. 
The size of the Sammy program is, up to a constant, the same as the number of rules in the Turing machine. That is 
$K_{Sammy}(F_s)= O(K_{Classical}(s))$ where $F_s$ is a functor that describes a string.  In a sense, this says that our Kolmogorov complexity is a generalization of classical 
Kolmogorov complexity.

We do not see why there should be a theorem that goes the other way. In other words, we do not think that 
a Turing machine can mimic an arbitrary Sammy program. If, in fact there are some categorical constructions that can be constructed by a Sammy program, but cannot be constructed by a Turing machine, then our Kolmogorov complexity is stronger than  classical Kolmogorov complexity theory.
Here is an example of a category and  a functor that can NOT be constructed by a Turing machine but might be able to be constructed by a Sammy program. Let ${\bf Halt}$ be a the ``halting category'' whose objects are the natural numbers and whose morphisms are defined below. Similarly there is the ``halting functor', $H$, from $\N$, the totally ordered category of the natural numbers, to $\overline{\bf{2}}$,  the category with two objects and a unique isomorphism between them, is defined on the right. 
$$ Hom_{\bf Halt}(n,n)= \left \{
\begin{array}{r@{\quad : \quad}l}
\mathbf \omega & \mbox{if } \varphi_n(n)\downarrow\\
Id_n & \mbox{if } \varphi_n(n)\uparrow \end{array}\right. \qquad \qquad  \qquad H(n)= \left \{
\begin{array}{r@{\quad : \quad}l}
1 & \mbox{if } \varphi_n(n)\downarrow\\
0 & \mbox{if } \varphi_n(n)\uparrow \end{array}\right.$$
Although, at present time, I do not know how to write a Sammy program to make such constructions, I believe that using infinite limits and colimits one should be able to build a type of infinite-time Turing machine to tell if regular Turing machines will halt or not. (However we are hesitant about making any conjectures. There is an interesting information-theoretic proof of the undecidability of the halting problem given on page 362 of  \cite{Calude}. Much work remains.)

Although we suspect that Sammy can actually program a larger class of functions than a Turing machine, however, there are some categorical constructions that are not programmable by Sammy (or any other language.) It is known that $K_{Classical}: Strings \lra \N$ is not a computable function. What about $K_{Sammy}$? First let us be careful about the definition of $K_{Sammy}$. It is a function that assigns to every category, functor, and natural transformation a natural number. We might as well assume that it only assigns natural transformations since identity natural transformations are simply functors and identity functors are simply categories. Let us think of $\Cat$ as the discrete category of natural transformation. We are going to forget the (two) composition structures on $\Cat$ because $K_{Sammy}$ does not behave well in terms of composition. So we have a functor $K_{Sammy}:\Cat \lra \N$. We prove that this functor is not constructible. The proof is a self-reference argument similar to the Berry paradox.
\begin{teo}
$K_{Sammy}:\Cat \lra \N$ is not constructible.
\end{teo}  

\noindent {\bf Proof.} Assume (wrongly) that $K=K_{Sammy}$ is, in fact, constructible, then there is a shortest program that describes $K$. In that case we can ask for the value of $K(K)$ (this is the core of self reference!). Let $K(K)=c$. Also, let $n$ be a natural number and let $P_n:\bf 1 \lra \N$ be a functor such that $P_n(0)=n$. Now use $K$ and and $P_n$ to construct the following pullback:
$$\xymatrix{\Cat_n \ar@{^{(}->}[r]\ar[d]& \Cat\ar[d]^K 
\\
 P_n\downarrow \N \ar@{^{(}->}[r]&\ N .
}$$
$ P_n \downarrow\N$ is the sub-total order of natural numbers that start at $n$. $\Cat_n$ is the discrete set of natural transformations whose shortest program is greater than or equal to $n$ operations. This pullback only needed a few more operations than $c$. Say that 
$K(\Cat_n| n)=c'$. However we can ``hardwire'' any $n$ into the program. If we do that, we get $K(\Cat_n)= c' + log~ n$.  
Choose an $n$ such that $n >> c' + log~ n$. Then $\Cat_n$ contains objects that require $n$ or more lines of code while we just described $\Cat_n$ in $c' + log~ n$ lines of code. This is like a Berry sentence. Contradiction! The only thing assumed is that $K$ was constructible. It is not constructible. QED.

We see this paper as just the beginning of a larger project to understand the complexity of categorical structures. Our work is far from done. With this notion of Kolmogorov complexity we get different notions of randomness, compressibility, and different notions of information. We would like to find upper bounds on some given categorical structures. We also would like to better clarify what is constructible and what is not. Another goal is to continue finding different categorical versions of the incompleteness theorems. We also would like to study different complexity measures. Rather than asking what is the shortest program that produces a categorical structure, we can ask how much time/space does a program take to create a certain structure. That is, what is the computational complexity of a structure. We can ask how much time does it take for the shortest program to produce that structure (logical depth.) All these measures induce hierarchies and classifications of categorical structures. There are also many other areas that we plan on studying. Here are a few.

There is a relationship between classical Kolmogorov complexity and Shannon's complexity theory. We would like to formulate a notion of Shannon's complexity theory for categories. There should be a definition of entropy of a category which should measure how rigid or flexible categorical structure is. Let 
$\C$ be a category, then $Aut(\C)$ is the group of automorphism functors $F:\C \lra \C$.
Define the ``entropy'' (or ``Hartley entropy'') of $\C$ as 
$H(\C)= Log_2 |Aut(\C)|$.
Just as there is a relationship between these measures for strings, there should be a relationship for categorical structures.

 So far we have restricted to classical categories, functors, and natural transformations. What about categories with more structure? For example, what can we say about a category that we know has all limits and colimits? What about enriched categories, higher categories, categories with structure, quasi-categories, etc? These different structures have been applied in almost every area of mathematics, computer science and theoretical physics. What we worked out above is only the first step. Such a study would be extremely interesting to shed some light on coherence theory. In this paper we saw that a pivotal fact of the Kolmogorov complexity of categories is that some categories are defined up to a unique isomorphism. Coherence theory generalizes such notions and is, in a sense, a higher dimensional version of uniqueness We will learn much about categorical information content and coherence theory by seeing the way they interact. 

This work should also be related to the important work in quantum information theory. We would like to study some of the physical and mathematical structures that occur in quantum mechanics with the developed Kolmogorov complexity tools.

Another area that we would like to explore is Occam’s razor \cite{OLR}. This is usually seen as a criteria in which to judge different physical theories. In short, physicists formulate functors
$F : $``Physical Phenomena'' $\lra$  ``Mathematical Structure.''
Universality of the theory demands that ``Physical Phenomena'' be as large as possible. In contrast, Occam's razor demands that ``Mathematical Structure" have low informational
content. We would like to use Kolmogorov complexity on both of these types of categories and the functors that relates them. We feel that with a better understanding of this we would be able to understand the question of why it seems that Occam's razor works so well.


\begin{thebibliography}{99}

\bibitem{Calude}Calude, Cristian, {\it Information and Randomness: An Algorithmic Perspective
Second Edition}
Springer-Verlag New York, 2002.
\bibitem{Hagino} Hagino, Tatsuya, {\sl A Categorical Programming Language}, available at \newline http://voxoz.com/publications/cat/Category
\bibitem{RandB}
Rydeheard, D.E., and Burstall, R.M., {\sl Computational Category Theory} available at http://www.cs.man.ac.uk/~david/categories/book/book.pdf 
\bibitem{Vitanyi}Li, Ming and     Vit\'anyi, Paul M. B.
{\sl An Introduction to Kolmogorov Complexity and its Applications. Second Edition}. Springer, 1997.  


\bibitem{OLR} Yanofsky, N.S., {\sl The Outer Limits of Reason:
What Science, Mathematics, and Logic Cannot Tell Us.} MIT Press, 2013. 
\bibitem{YanoFut} Yanofsky, N.S., ``Algorithmic Information Theory in Categorical Algebra'' work in progress. 
\end{thebibliography}
\end{document}